\documentclass{article}
\usepackage[utf8]{inputenc}
\usepackage{authblk}
\usepackage[english]{babel}
\usepackage{amsmath}
\usepackage{amssymb}
\usepackage{amsfonts}
\usepackage{amsthm}
\usepackage{hyperref}
\usepackage{breqn}

\newtheorem{theorem}{Theorem}[section]

\def\ff{\ensuremath{\mathbf{F}}}

\title{Collection of Polynomials over  Finite Fields Providing Involutary Permutations}

\author{
  Kevinsam B and  P Vanchinathan\\
  School of Advanced Sciences\\
  Vellore Institute of Technology\\
  Chennai 600 127\\ INDIA\\
  \texttt{kevinsam.b2020@vitstudent.ac.in\qquad
	  vanchinathan.p@vit.ac.in}
}

\begin{document}
\maketitle

\begin{abstract} 
  \noindent For an odd prime power $q$ satisfying $q\equiv 1\pmod 3$ we construct
	totally $2(q-1) $ permutation polyomials, all giving involutory
	permutations with exactly $ 1+ \frac{q-1}3$ fixed points. Among them $(q-1)$ 
	polynomials are trinomials, and the rest are 6-term polynomials.

\end{abstract}

\section{Introduction}
An immediate consequence of Lagrange interpolation theorem  is that  every function from a finite 
field $\mathbf{F}_q$  to itself is a polynomial function of degree  at most $q-1$.
Then it is natural to ask the question which  of those  polynomials will represent 
bijective functions.
Such polynomials named Permutation polynomials have found many applications
in cryptography, coding theory and  combinatorics.
Lidl and Niederreiter's book \cite{lidl} in Chapter 7 deals with permutation polynomials.
Chapter 8 of \cite{handbook} by Mullen and Panerio has extensive results with large bibliography.
The survey article by Hou \cite{hou} is more recent and gives the staus upto 2015 on pemutation polynomials.

One natural problem here is given a  permutation  with specific
cycle decomposition (which means specifying a conjugacy class in the symmetric
group of all permutations of a finite field)  how the polynomials
corresponding to them can be identified. For example, if there is 
 any pattern in their coefficients.  
Simplest kind of permutations, namely transpositions,  as polynomials, have 
degree $q-2$ (the maximum for a permutation polynomial) 
and often have as many terms present in general.
However for cyclic permutations of length $q$,   the polynomials representing them can be 
as low as degree 1, or  very high depending on $q$. The relationship between the cycle
type of a permutation and the kind of polynomial representing it seems tricky to pinpoint.
This relationship is understood only in very limited cases.

\textit{In general constructing a family of polynomials all yielding permutations with  same  
cycle decomposition is a challenging task.}  (See the introduction in \cite{redei}). 

The second author of this paper jointly with Anitha G has in a recent arXiv preprint
constructed many families of binomial and trinomial permutations with cycle structure as $1+m^d$ and more.
See \cite{anitha}.

In this paper we have succeeded in addressing this task for certain involutions.
Another desirable property is that our polynomials are all trinomials, or 6-term polynomials.

More precisely, for fields $\ff_q$  with $q$ odd and $q\equiv1\pmod3$, we have constructed 
$2(q-1)$ polynomials 
with cycle structure $1^{(q+2)/3} + 2^{(q-1)/3 }$
half of them with three terms (trinomials)  and remaining with six terms.
(for precise statements see Theorems 2.1--2.3  for trinomials and  Theorems 2.4--2.6 for six-term polynomials).

Many authors have constructed 
 permutation trinomials where the focus is not
 on the cycle type of the permutations. These constructions  can be found in
 the works of Bartoli and Timponella \cite{bart}, 
 Gupta and Sharma\cite{rohit} (in characteristic 2),
Hou, \cite{hou1,hou2},  Lavorante\cite{lavo}, and Li, Qu and Chen \cite{lqc}.

\vspace{1pc}
NOTATIONS:

$q$ an odd prime power, always assumed to be   ${}\equiv1\pmod 3.$

$m = (q-1)/3$.

$\ff_q$ a finite field of order $q$.

$\gamma$ a generator of $\ff_q^*$.

\section{Main results}
With the notation as in the introduction we are ready state our results.
\begin{theorem}\label{th1}
For  an integer $k$ with $0\leq k\leq m-1$,  
	 define $a_k, b_k,c_k \in \mathbf{F}_{q}^*$ as below:

	 \vskip3pt
$\displaystyle a_k =\frac{\gamma^{2m+6k+2}+\gamma^{m+3k+1}+1}{3\gamma^{m+3k+1}} $

$\displaystyle b_k =\frac{\gamma^{2m}+\gamma^{m+3k+1}+\gamma^{6k+2}}{3\gamma^{m+3k+1}}$ ~and  

$\displaystyle c_k =\frac{\gamma^{6k+2}+\gamma^{3k+1}+1}{3\gamma^{3k+1}}$

\noindent
The  trinomials with above $a_k,b_k,c_k$ as coefficients, defined as 
$$g_k(x) = a_kx^{2m+1} + b_kx^{m+1} + c_kx $$ 
are all  permutation polynomials over $\mathbf{F}_{q}$,
and  they all give involutory permutation of $\mathbf{F}_{q}$ with exactly $m+1$ fixed points.
\end{theorem}

\noindent \textit{Remark:} To avoid clutter  in the notation, hereafterwards
we drop the subscript $k$ and simply use  $a,b,c$ and $g(x)$.
So this theorem constructs $m = (q-1)/3$ permutation trinomials  representing
involutions of $\mathbf{F}_q$ with exactly $(q+2)/3$ fixed-points.

\begin{proof} 
	We will prove  $g(x)$ is an  involution of $\mathbf{F}_{q}$
	by explicitly describing all its fixed points 
	and computing its values on other elements.
More precisely we show, for $i=0,1,2,\ldots, m-1$, that
\begin{itemize}
	\item[(i)] $g(\gamma^{3i})   = \gamma^{3i}$
	\item[(ii)] $g(\gamma^{3i+1}) = \gamma^{3(i+k)+2}$ and
	\item[(iii)] $g(\gamma^{3(i+k)+2}) = \gamma^{3i+1}$
\end{itemize}

\noindent \textit{Proof of (i)}:
\begin{flalign*}
g(\gamma^{3i}) & = a(\gamma^{3i})^{2m+1} + b(\gamma^{3i})^{m+1} + c(\gamma^{3i}) & \\
 & =\gamma^{3i}\left[ a(\gamma^{3i})^{2m} + b(\gamma^{3i})^m + c \right] & \\
& =\gamma^{3i}\left[ a + b + c \right] (\mbox{ since } m= (q-1)/3 \mbox{ and } \gamma^{3m}=1). 
 \end{flalign*}
	Now proof of (i) will be complete if we show $a+b+c=1$.
	First  we observe  that $\gamma^m$ being an element of order 3, satisfies
	$\gamma^{2m} + \gamma^m + 1=0$. Now we have, by definition,
 \begin{flalign*}
a+b+c & = \frac{\gamma^{2m+6k+2}+\gamma^{m+3k+1}+1}{3\gamma^{m+3k+1}} + 
	   \frac{\gamma^{2m}+\gamma^{m+3k+1}+\gamma^{6k+2}}{3\gamma^{m+3k+1}} \\
	 &    \hspace{12.5em}+ \frac{\gamma^{6k+2}+\gamma^{3k+1}+1}{3\gamma^{3k+1}} & \\
	  & = \frac{3\gamma^{m+3k+1}+(\gamma^{6k+2}+1)(\gamma^{2m} + \gamma^m + 1)}{3\gamma^{m+3k+1}} & \\
& =\frac{3\gamma^{m+3k+1}}{3\gamma^{m+3k+1}}   = 1.
 \end{flalign*}
This shows that the $g(x)$ has at least $m+1$ fixed points, along with 0.\\

\noindent \textit{Proof of (ii)}:
	To show  $g(\gamma^{3i+1}) = \gamma^{3(i+k)+2}$ we proceed to compute
	as below:
\begin{flalign*}
g(\gamma^{3i+1}) &= a(\gamma^{3i+1})^{2m+1} + b(\gamma^{3i+1})^{m+1} + c(\gamma^{3i+1}) & \\
 &=\gamma^{3i+1}\left[ a(\gamma^{3i+1})^{2m} + b(\gamma^{3i+1})^m + c \right] & \\
 &=\gamma^{3i+1}\left[ a\gamma^{2m} + b\gamma^m + c \right]
\end{flalign*}
Substituting the values of $a,b$ and $c$ this becomes
\begin{dmath*}
=\gamma^{3i+1} \left[ \frac{\gamma^{2m+6k+2}+\gamma^{m+3k+1}+1}{3\gamma^{m+3k+1}}(\gamma^{2m}) + \frac{\gamma^{2m}+\gamma^{m+3k+1}+\gamma^{6k+2}}{3\gamma^{m+3k+1}}(\gamma^{m})\\ + \frac{\gamma^{6k+2}+\gamma^{3k+1}+1}{3\gamma^{3k+1}}\right] \\
=\gamma^{3i+1} \left[ \frac{\gamma^{2m+6k+2}+\gamma^{m+3k+1}+1}{3\gamma^{3k+1}}(\gamma^{m}) + \frac{\gamma^{2m}+\gamma^{m+3k+1}+\gamma^{6k+2}}{3\gamma^{3k+1}}\\ + \frac{\gamma^{6k+2}+\gamma^{3k+1}+1}{3\gamma^{3k+1}}\right] \\
=\frac{\gamma^{3i+1}}{3\gamma^{3k+1}} \left[\gamma^{3m+6k+2}+\gamma^{2m+3k+1}+\gamma^{m} + \gamma^{2m}+\gamma^{m+3k+1}+\gamma^{6k+2} + \gamma^{6k+2}+\gamma^{3k+1}+1\right] \\
=\frac{\gamma^{3i+1}}{3\gamma^{3k+1}} \left[ 3\gamma^{6k+2} + (\gamma^{3k+1} + 1)(\gamma^{2m} + \gamma^m +1)\right]
\end{dmath*}
Again using the fact that $\gamma^{2m}+\gamma^{m} + 1 =0$,
 the last line above simplifies to 
\begin{flalign*}
&\hspace{1em}=\frac{\gamma^{3i+1}}{3\gamma^{3k+1}} \left[ 3\gamma^{6k+2}  \right] &\\
&\hspace{1em}=\gamma^{3(i+k)+2}.
\end{flalign*}
\noindent \textit{Proof of (iii)}:
	To show $g(\gamma^{3(i+k)+2})= \gamma^{3k+1}$ we proceed
	as below:
 \begin{flalign*}
g(\gamma^{3(i+k)+2}) & = a(\gamma^{3(i+k)+2})^{2m+1} + b(\gamma^{3(i+k)+2})^{m+1} + c(\gamma^{3(i+k)+2})& \\
&=\gamma^{3(i+k)+2}\left[ a(\gamma^{3k+2})^{2m} + b(\gamma^{3k+2})^m + c \right]& \\
&=\gamma^{3(i+k)+2}\left[ a\gamma^{m} + b\gamma^{2m} + c \right]
 \end{flalign*}
Substituting the values of $a,b$ and $c$ this becomes

\begin{dmath*}
=\gamma^{3(i+k)+2} \left[ \frac{\gamma^{2m+6k+2}+\gamma^{m+3k+1}+1}{3\gamma^{m+3k+1}}(\gamma^{m}) + \frac{\gamma^{2m}+\gamma^{m+3k+1}+\gamma^{6k+2}}{3\gamma^{m+3k+1}}(\gamma^{2m})\\ + \frac{\gamma^{6k+2}+\gamma^{3k+1}+1}{3\gamma^{3k+1}}\right] \\
=\gamma^{3(i+k)+2} \left[ \frac{\gamma^{2m+6k+2}+\gamma^{m+3k+1}+1}{3\gamma^{3k+1}} + \frac{\gamma^{2m}+\gamma^{m+3k+1}+\gamma^{6k+2}}{3\gamma^{3k+1}}(\gamma^{m})\\ + \frac{\gamma^{6k+2}+\gamma^{3k+1}+1}{3\gamma^{3k+1}}\right] \\
=\frac{\gamma^{3(i+k)+2}}{3\gamma^{3k+1}} \left[\gamma^{2m+6k+2}+\gamma^{m+3k+1}+1 + \gamma^{3m}+\gamma^{2m+3k+1}+\gamma^{m+6k+2}\\ + \gamma^{6k+2}+\gamma^{3k+1}+1\right] \\
=\frac{\gamma^{3(i+k)+2}}{3\gamma^{3k+1}} \left[ 3 + (\gamma^{6k+2} + \gamma^{3k+1})(\gamma^{2m} + \gamma^m +1)\right]
\end{dmath*}
 Again using the fact that $\gamma^{2m}+\gamma^{m} + 1 =0$ the above becomes
 \begin{flalign*}
&\hspace{1em}=\frac{\gamma^{3(i+k)+2}}{3\gamma^{3k+1}} \cdot  3 & \\
&\hspace{1em}=\gamma^{3i+1} 
 \end{flalign*}
\end{proof}
\noindent\textbf{Remark:}\quad
The above permutation polynomial is a cyclotomic mapping fixing pointwise
a subgroup $H$ of index 3 in $\mathbf{F}_q^*$ and swapping elements
of the coset $\gamma H$ with ``corresponding''  elements of the coset $\gamma^2H$.

\textit{The word corresponding given in quotes above  means `an element
 k steps forward' in the other coset}.

In general we can keep any one of the three cosets of $H$ pointwise fixed and
swap the corresponding elements of the other two cosets
again getting  involutions with exactly $m+1$ fixed points.

We were also able to find the permutation polynomial for those involutions too.

They also, fortunately, turned out to be  trinomials involving the same terms  
but with different choices for the coefficients $a,b,c$. The next two theorems 
describe those trinomials. 

\begin{theorem}
	For integers $k$ with $0\leq k\leq m-1$ define $a, b,c \in \mathbf{F}_{q}^*$ as below: \\
$\displaystyle  a=\frac{\gamma^{2m+3k+2}+\gamma^{m+6k+4}+1}{3\gamma^{m+3k+2}}$ \\
$\displaystyle  b=\frac{\gamma^{2m}+\gamma^{m+6k+4}+\gamma^{3k+2}}{3\gamma^{m+3k+2}}$ ~and \\
$\displaystyle  c=\frac{\gamma^{6k+4}+\gamma^{3k+2}+1}{3\gamma^{3k+2}}$ \\
	Then the polynomial $g(x) = ax^{2m+1} + bx^{m+1} + cx $ is a
 permutation polynomial over $\mathbf{F}_{q}$. Moreover, $g(x)$ represents an involutory permutation of $\mathbf{F}_{q}$ with exactly $m+1$ fixed points.
\end{theorem}

\begin{proof}
	We claim that
\begin{itemize}
	\item[(i)] $g(\gamma^{3i+1}) = \gamma^{3i+1}$
	\item[(ii)] $g(\gamma^{3i}) = \gamma^{3(i+k)+2}$ and
	\item[(iii)] $g(\gamma^{3(i+k)+2}) = \gamma^{3i}$
\end{itemize}
where $i=0,1,2,\dots, m-1$. The proof uses the same kind of arguments as in Theorem~\ref{th1}
	and so omitted. 
\end{proof}

\begin{theorem}
	Define $a, b,c \in \mathbf{F}_{q}^*$ as below: \\
$\displaystyle  a=\frac{\gamma^{2m}+\gamma^{m+6k+2}+\gamma^{3k+1}}{3\gamma^{m+3k+1}}$ \\
$\displaystyle  b=\frac{\gamma^{2m+3k+1}+\gamma^{m+6k+2}+1}{3\gamma^{m+3k+1}}$ ~and \\
$\displaystyle  c=\frac{\gamma^{6k+2}+\gamma^{3k+1}+1}{3\gamma^{3k+1}}$ \\
for each $k=0,1,\dots,m-1$. Then the polynomial $g(x) = ax^{2m+1} + bx^{m+1} + cx $ is a
	permutation polynomial over $\mathbf{F}_{q}$, and  $g(x)$ gives  an 
	involutory permutation of $\mathbf{F}_{q}$ with exactly $m+1$ fixed points.
\end{theorem}
\begin{proof}
In fact more specifically
	\begin{itemize}
		\item[(i)] $g(\gamma^{3i+2}) = \gamma^{3i+2}$
		\item[(ii)] $g(\gamma^{3i}) = \gamma^{3(i+k)+1}$ and
		\item[(iii)] $g(\gamma^{3(i+k)+1}) = \gamma^{3i}$
\end{itemize}
where $i=0,1,2,\dots, m-1$.
	Once again we omit the proof as it employs  arguments similar to the ones 
	presented in the proof of Theorem~\ref{th1}.\\
\end{proof}


\noindent\textbf{Remark:}\quad  Combining the above three theorems we have
produced $q-1$ permutation trinomials  represnting 
involutions, each  with $(q+2)/3$ fixed points 
for fields of order $q\equiv1\pmod 3$.

Now we move onto  permutation polynomials of 6 terms that is a cyclotomic
mapping for the subgroup $H$ mentioned above, which is also an involution.

\begin{theorem}\label{th2}
	Define $a, b,c \in \mathbf{F}_{q}^\ast$ as below: \\
$$\displaystyle  a=\frac{2\gamma^{3k}}{3}, \quad  b=\frac{1}{3}\quad \text{and}\quad c=\frac{-\gamma^{3k}}{3}$$
where $k=0,1,\dots ,m-1$. Then the polynomial
$$g(x) = ax^{3m-1}+bx^{2m+1}+cx^{2m-1}  + bx^{m+1}  +c x^{m-1} +bx$$

	is permutation polynomials over $\mathbf{F}_{q}$, actually $g(x)$ represents an involutory permutation of $\mathbf{F}_{q}$ with exactly $m+1$ fixed points.

\begin{proof}
	First note that $g(x)$ has zero as one fixed point. 
To prove that $g(x)$ is involutory permutation of $\mathbf{F}_{q}$,
	 we will show
\begin{itemize}
	\item[(i)] $g(\gamma^{3i})= \gamma^{3i}$
	\item[(ii)] $g(\gamma^{3i+1})= \gamma^{3(k-i)-1}$ and
	\item[(iii)] $g(\gamma^{3(k-i)-1})= \gamma^{3i+1}$
\end{itemize}
where $i=0,1,2,\dots, m-1$\\
The polynomial can be written as
$$g(x) = ax^{3m-1}+b\left[x^{2m+1} + x^{m+1} + x\right] + c\left[x^{2m-1} + x^{m-1}\right]$$
	Proof of (i): $g(\gamma^{3i})\mapsto \gamma^{3i}$\\
\begin{flalign*}
g(\gamma^{3i}) &= a(\gamma^{3i})^{3m-1}+b\left[(\gamma^{3i})^{2m+1} + (\gamma^{3i})^{m+1} + \gamma^{3i}\right]& \\
		&\hphantom{abcdekdbqkidd} + c\left[(\gamma^{3i})^{2m-1} + (\gamma^{3i})^{m-1}\right]& \\
	&=a(\gamma^{3i})^{-1}+b\gamma^{3i}\left[(\gamma^{3i})^{2m}+(\gamma^{3i})^{m}+1\right]+c\gamma^{-3i}\left[(\gamma^{3i})^{2m}+(\gamma^{3i})^{m}\right]& \\
&=a(\gamma^{-3i})+b\gamma^{3i}\left[1+1+1 \right] + c\gamma^{-3i}\left[ 1+1 \right]
\end{flalign*}
Now substituting the values for the coefficients $a,b,c$ we get
\begin{flalign*}
&\hspace{2.5em}=\frac{2\gamma^{3k}}{3}(\gamma^{-3i}) +  \frac{1}{3}(\gamma^{3i})(3) - \frac{\gamma^{3k}}{3}(\gamma^{-3i})(2)& \\
&\hspace{2.5em}= \frac{2}{3}\gamma^{3(k-i)} + \gamma^{3i} - \frac{2}{3}\gamma^{3(k-i)}& \\
&\hspace{2.5em}= \gamma^{3i}
\end{flalign*}
This shows that $g(x)$ has $m+1$ at least fixed points. That they are the only fixed
points will follow from the remaining part of proof.

	Proof of (ii): $g(\gamma^{3i+1})= \gamma^{3(k-i)-1}$\\
 \begin{flalign*}
g(\gamma^{3i+1}) &= a(\gamma^{3i+1})^{3m-1}+b\left[(\gamma^{3i+1})^{2m+1} + (\gamma^{3i+1})^{m+1} + \gamma^{3i+1}\right] & \\
		&\hphantom{abcdekdbqkidddd}+ c\left[(\gamma^{3i+1})^{2m-1} + (\gamma^{3i+1})^{m-1}\right]& \\
& = a(\gamma^{3k+1})^{-1}+b(\gamma^{3i+1})\left[(\gamma^{3i+1})^{2m} + (\gamma^{3i+1})^{m} + 1\right]& \\
		&\hphantom{abcdekdbqknis} + c(\gamma^{-3i-1})\left[(\gamma^{3i+1})^{2m} + (\gamma^{3i+1})^{m}\right]& \\
& = a(\gamma^{-3i-1})+ b(\gamma^{3i+1})\left[\gamma^{2m} + \gamma^{m} + 1\right] + c (\gamma^{-3i-1})\left[\gamma^{2m} + \gamma^{m}\right]
 \end{flalign*}
Since $\gamma^m$ be a primitive third root of unity, we have $\gamma^{2m}+\gamma^{m} + 1 =0$ and then substituting the values of $a,b$ and $c$ we get
\begin{flalign*}
& \hspace{3.5em} = \frac{2\gamma^{3k}}{3}(\gamma^{-3i-1}) - \frac{\gamma^{3k}}{3}(\gamma^{-3i-1})\left(-1 \right)& \\
& \hspace{3.5em}= \frac{3\gamma^{3(k-i)-1}}{3}& \\
& \hspace{3.5em}= \gamma^{3(k-i)-1}
\end{flalign*}
	\textit{Proof of (iii)}: $g(\gamma^{3(k-i)-1}) = \gamma^{3i+1}$
	\begin{flalign*}
		g(\gamma^{3(k-i)-1})&= a(\gamma^{-3(k-i)+1})+b\gamma^{3(k-i)-1}\left[(\gamma^{3(k-i)-1})^{2m} + (\gamma^{3(k-i)-1})^{m} + 1\right]& \\ 
		&\hspace{7.5em}+ c\gamma^{-3(k-i)+1}\left[(\gamma^{3(k-i)-1})^{2m} + (\gamma^{3(k-i)-1})^{m}\right] & \\
		&= a\gamma^{-3(k-i)+1}+b\gamma^{3(k-i)-1}\left[\gamma^{-2m} + \gamma^{-m} + 1\right]& \\ 
		&\hspace{6.5em} + c\gamma^{-3(k-i)+1}\left[\gamma^{-2m} + \gamma^{-m}\right]& \\
		&= a\gamma^{-3(k-i)+1}+b\gamma^{3(k-i)-1}\left[\gamma^{m} + \gamma^{2m} + 1\right]& \\ 
		&\hspace{6.5em} + c\gamma^{-3(k-i)+1} \left[\gamma^{m} + \gamma^{2m}\right]
\end{flalign*}
Since $\gamma^m$ be a primitive third root of unity, we have $\gamma^{2m}+\gamma^{m} + 1 =0$ and then substituting the values of $a,b$ and $c$ this becomes
	\begin{flalign*}
		&\hspace{5.4em}= \frac{2\gamma^{3k}}{3}\gamma^{-3(k-i)+1} - \frac{\gamma^{3k}}{3} \gamma^{-3(k-i)+1}\left[-1\right]& \\
		&\hspace{5.4em}= \frac{2\gamma^{3i+1}}{3} + \frac{\gamma^{3i+1}}{3}& \\
		&\hspace{5.4em}= \gamma^{3i+1}.
	\end{flalign*}
\end{proof}
\end{theorem}


\begin{theorem}
	Define $a, b,c,d,e,f \in \mathbf{F}_{q}^\ast$ as below: 

\vspace{5pt}
\noindent$\displaystyle a=\frac{2\gamma^{3k+2}}{3},    \hspace{1.34pt}
b=\frac{\gamma^{m}}{3},\hspace{1.34pt}
c=\frac{-\gamma^{m+3k+2}}{3},\hspace{1.34pt}
d=\frac{\gamma^{2m}}{3},\hspace{1.34pt}
e=\frac{-\gamma^{2m+3k+2}}{3}~\text{and} 
	~ f=\frac{1}{3}$\\[5pt]
where $k=0,1,\dots ,m-1$. Then the polynomial $$g(x) = ax^{3m-1}+bx^{2m+1}+cx^{2m-1}+ dx^{m+1} + ex^{m-1} + fx$$ is a permutation polynomial over $\mathbf{F}_{q}$, actually $g(x)$ represents an involutory permutation of $\mathbf{F}_{q}$ with exactly $m+1$ fixed points.

\end{theorem}
\begin{proof}
To prove that $g(x)$ is involutory permutation of $\mathbf{F}_{q}$. $g(x)$ has zero as one fixed point. Now we need to prove the following three cases:
\begin{itemize}
	\item[(i)] $g(\gamma^{3i+1}) = \gamma^{3i+1}$
	\item[(ii)] $g(\gamma^{3i}) =  \gamma^{3(k-i)+2}$ and
	\item[(iii)] $g(\gamma^{3(k-i)+2}) = \gamma^{3i}$
\end{itemize}
where $i=0,1,2,\dots, m-1$. The proof uses the same kind of arguments as in Theorem~\ref{th2}.\\
\end{proof}


\begin{theorem}
	Define $a, b,c,d,e,f \in \mathbf{F}_{q}^\ast$ as below: \\
$$\displaystyle a=\frac{2\gamma^{3k+1}}{3},~
b=\frac{\gamma^{2m}}{3},~
c=\frac{-\gamma^{2m+3k+1}}{3},~
d=\frac{\gamma^{m}}{3},~
e=\frac{-\gamma^{m+3k+1}}{3}~ and ~
f=\frac{1}{3}$$\\
where $k=0,1,\dots ,m-1$. Then the polynomial $$g(x) = ax^{3m-1}+bx^{2m+1}+cx^{2m-1}+ dx^{m+1} + ex^{m-1} + fx$$ is permutation polynomials over $\mathbf{F}_{q}$, actually $g(x)$ represents an involutory permutation of $\mathbf{F}_{q}$ with exactly $m+1$ fixed points.

\end{theorem}
\begin{proof}
To prove that $g(x)$ is involutory permutation of $\mathbf{F}_{q}$. $g(x)$ has zero as one fixed point. Now we need to prove the following three cases:
\begin{itemize}
	\item[(i)] $g(\gamma^{3i+2}) =  \gamma^{3i+2}$
	\item[(ii)] $g(\gamma^{3i}) =  \gamma^{3(k-i)+1}$ and
	\item[(iii)] $g(\gamma^{3(k-i)+1}) =  \gamma^{3i}$
\end{itemize}
where $i=0,1,2,\dots, m-1$. The proof uses the same kind of arguments as in Theorem~\ref{th2}\\
\end{proof}

\section{Conclusion}
To summarise, for a field of $q$ elements with $q$ odd and $q\equiv1\pmod3,$
we have constructed $2(q-1)$ permutation polynomials
all representing permutations in a single conjugacy class
given by the integer partition
\[q=\underbrace{1+1+\cdots+1}_{(q+2)/3\ \rm terms} +\underbrace{2+2+
\cdots+2}_{(q-1)/3\ \rm terms}\]
It is to be noted that the polynomials we have constructed have  all  been described 
with respect to  a specific choice
of generator $\gamma$ of $\ff_q^*$. Varying this choice might lead to bigger
collection of permutation polynomials with the same cycle decomposition.

That these collections are disjoint for different choices of $\gamma$ remains to be
investigated.

\end{document}